\documentclass[11pt,oneside]{amsart}
\usepackage[latin9]{inputenc}
\usepackage{geometry}
\usepackage{textcomp}
\usepackage{amsthm}
\usepackage{amstext}
\usepackage{amssymb}

\makeatletter

\newcommand{\lyxmathsym}[1]{\ifmmode\begingroup\def\b@ld{bold}
  \text{\ifx\math@version\b@ld\bfseries\fi#1}\endgroup\else#1\fi}

\providecommand{\tabularnewline}{\\}

\numberwithin{equation}{section}
\numberwithin{figure}{section}
\theoremstyle{plain}
\newtheorem{thm}{\protect\theoremname}
  \theoremstyle{definition}
  \newtheorem{example}[thm]{\protect\examplename}

\title{}

\usepackage{amsfonts}\usepackage{enumerate}

\newtheorem{theorem}{Theorem}
\theoremstyle{plain}

\newtheorem{corollary}{Corollary}

\newtheorem{definition}{Definition}

\numberwithin{equation}{section}

\makeatother

  \providecommand{\examplename}{Example}
\providecommand{\theoremname}{Theorem}

\begin{document}

\title{On Modulo AG-groupoids}

\author{Amanullah$^{1}$}

\email{amanswt@hotmail.com}

\author{M. Rashad$^{1}$}

\email{rashad@uom.edu.pk}

\author{I. Ahmad$^{1,*}$}

\email{iahmaad@hotmail.com}

\author{M. Shah$^{2}$}

\email{shahmaths\_problem@hotmail.com}

\address{1. Department of Mathematics University of Malakand, Chakdara, Pakistan.}

\address{2. Department of Mathematics, Govt Post Graduate College Mardan,
Pakistan.}

\keywords{AG-groupoids $(\mbox{mod }n)$; AG-groups $(\mbox{mod }n)$; construction;
$T^{3}$-AG-groupoid; cancellative AG-groupoid.\\
{*}Corresponding author}
\begin{abstract}
A groupoid $G$ is called an AG-groupoid if it satisfies the left
invertive law: $(ab)c=(cb)a$. An \textit{\emph{AG-group}} $G$,\textit{\emph{
is an AG-groupoid}} with left identity\emph{ }$e\in G$\emph{ }(that
is, $ea=a$\emph{ }for all\emph{ }$a\in G$) and for all $a\in G$
there exists $a^{^{-1}}\in G$ such that\emph{ }$a^{^{-1}}a=aa^{-1}=e$.
In this article we introduce the concept of AG-groupoids $(\mbox{mod }n)$
and AG-group $(\mbox{mod }n)$ using Vasantha's constructions \cite{WBV}.
This enables us to prove that AG-groupoids $(\mbox{mod }n)$ and AG-groups
$(\mbox{mod }n)$ exist for every integer $n\geq3$. We also give
some nice characterizations of some classes of AG-groupoids in terms
of AG-groupoids $(\mbox{mod }n)$.
\end{abstract}
\maketitle

\section{Introduction}

\noindent Construction for any algebraic structure is always very
important for its developement. The examples so obtained are sometimes
not even possible through computers to come by. Open problems and
conjunctures are often answered by constructing examples for them.
Several construction are available for forming quasigroup and loops.
For example, an infinite family of nonassociative noncommutative C-loops
whose smallest member is the smallest non-associative noncommutative
C-loop of order 12 has been constructed in \cite{PV1}. We can obtained
manually a C-loop of this family of much much higher order that might
not be possible through computers. Sometimes a construction can be
implemented in computer which then makes the job easier. For example,
a construction of AG-groups from abelian groups has been implemented
in GAP, the details of this construction and its implementation can
be found in \cite{AG-groups}. M. S. Kamran has also discussed another
construction of AG-groups from abelian groups in his PhD thesis \cite{SKamran}.
Several types of a structure can be obtained from each other through
some specific constructions. For example, this has been done for AG-groupoids
\cite{IMS}. Several constructions of groupoids have been done by
W. B. Vasantha \cite{WBV}. In this paper we extend Vasantha's constructions
to AG-groupoids by imposing some conditions on them. The structure
of AG-groupoid is considered one of the most interesting structure
among the non-associative structures. A considerable achivement has
been done for the improvement of AG-groupoids by various researchers
see for instance\cite{intro,key-1,RAAS,AGst,AGst-1}. This extension
will really give another push to the study of AG-groupoids which gets
broadened quite rapidly these days. This will also make Vasantha's
constructions of groupoids more valuable. The paper can be considered
as a sort of applications of the mentioned constructions. It also
gives some nice characterizations of some classes of AG-groupoids
in new scenario, that is, in terms of AG-groupoids $(\mbox{mod }n)$.

\noindent In Section 2 we will establish existence of AG-groupoid
$(\mbox{mod }n)$ and in Section 3 we will construct AG-groups $(\mbox{mod }n)$
and will provide some finite examples to show their existence.

\noindent An AG-groupoid (or LA-semigroup) $G$ is a groupoid in which
the left invertive law: $(ab)c=(cb)a$ holds \cite{KN}. An AG-groupoid
is a generalization of a commutative semigroup and an AG-group is
a generalization of abelian group. Recently some new classes of AG-groupoids
have been discovered in \cite{SIA} and \cite{MshahT}. However we
will need the following definitions. 

\noindent An \textit{\emph{AG-groupoid}} $G$ is called :

\noindent \begin{enumerate}[(i)]

\item an AG-band, if $a^{2}=a\,\,\forall a\in G$ \cite{SP};

\item an AG-group, if $b\ast a=c\ast a\Rightarrow a\ast b=a\ast c$;

\item a $T_{l}^{3}$-AG-groupoid, if $a\ast b=a\ast c\Rightarrow b\ast a=c\ast a$;

\item a $T_{r}^{3}$-AG-groupoid, if $b\ast a=c\ast a\Rightarrow a\ast b=a\ast c$;

\item a transitively commutative AG-groupoid, if $a\ast b=b\ast a\mbox{ and }b\ast c=c\ast b\Rightarrow a\ast c=c\ast a$;

\item a cancellative AG-groupoid, if $a\ast x=a\ast y\Rightarrow x=y$.

\noindent \end {enumerate}

\noindent The following definitions are needed from \cite{WBV}. 

\noindent \begin{definition} \label{Gpd1} Let $Z_{n}=\left\{ 0,1,2,\cdots,n\lyxmathsym{\textendash}1\right\} ;n\ge3$.
For $a,b\in Z_{n}\setminus\{0\},$ define a binary operation $\lyxmathsym{\textasteriskcentered}$
on $Z_{n}$ as follows $a*b=ta+ub(\mbox{mod }n)$ where $t,u$ are
two distinct elements in $Z_{n}\setminus\{0\}$ and $(t,u)=1$ here
` + ' is the usual addition of two integers and ` $ta$ ' means the
product of the two integers ` $t$ ' and ` $a$ '. This groupoid will
be denoted by $\left(Z_{n},\,(t,u),\,\lyxmathsym{\textasteriskcentered}\right)$
or in short by $Z_{n}(t,u)$. By varying $t,u\in Z_{n}\setminus\{0\}$
with $(t,u)=1$ we get a collection of groupoids for a fixed integer
$n$. This collection of groupoids is denoted by $Z(n)$ that is $Z(n)=\left\{ \left(Z_{n},(t,u),\lyxmathsym{\textasteriskcentered}\right)\,|\,\mbox{for distinct integers }t,u\in Z_{n}\setminus\{0\}\mbox{\,\ such that }(t,u)=1\right\} $.
Clearly every groupoid in this class is of order $n$. \end{definition}

\noindent \begin{definition} \label{Gpd2} If $(t,u)$ need not always
be relative prime but $t\neq u$ and $t,u\in Z_{n}\setminus\{0\}$
in Definition \ref{Gpd1} we get a new extended class of $Z(n)$ denoted
by $Z^{*}(n)$. \end{definition}

\noindent \begin{definition} \label{Gpd3} If $(t,u)$ need not always
be distinct in Definition \ref{Gpd1} we get a new enlarge class of
$Z^{*}(n)$ denoted by $Z^{**}(n)$. \end{definition} 

\noindent \begin{definition} \label{Gpd4} If $t,u\in Z_{n}$ where
$t$ or $u$ can also be zero in Definition \ref{Gpd1} we get a new
class denoted by $Z^{***}(n)$ contains $Z^{**}(n)$ . \end{definition}

\section{\noindent Existence of AG-groupoids $(\mbox{mod}\, n)$}

In this section, we introduce AG-groupoids $(\mbox{mod }n)$ as a
subclass of the class $Z^{***}(n)$. We study these AG-groupoids $(\mbox{mod }n)$
and obtain some results about them. The following theorem guarantees
the existence of AG-groupoids $(\mbox{mod }n)$ for $n\geq3$, indeed
it provides us with a simple way of construction of AG-groupoids $(\mbox{mod }n)$
of any finite order.

\noindent \begin {theorem}\label{GPDS} Let $Z_{n}=\{0,1,2,\cdots,n-1\},\, n\geq3,\, n<\infty$.
A groupoid in $Z^{***}(n)$ is an AG-groupoid if $t^{2}\cong u(\mbox{mod}\, n)\,\,\mbox{for any }t,u\in Z_{n}$.\end {theorem}

\noindent \begin {proof}Let $Z_{n}=\{0,1,2,\cdots,n-1\},\, n\geq3,\, n<\infty$;
satisfies $t^{2}\cong u(\mbox{mod}\, n)$ for any $t,u\in Z_{n}\setminus\{0\}$.
To show $Z_{n}(t,u)$ is an AG-groupoid, we have to show that the
left invertive law, that is, $(a\cdot b)\cdot c=(c\cdot b)\cdot a\,\,\forall\, a,b,c\in Z_{n}$
holds. Now
\begin{eqnarray*}
(a\cdot b)\cdot c & \cong & (t(ta+ub)+uc)(\mbox{mod}\, n)\\
 & \cong & (t^{2}a+tub+uc)(\mbox{mod}\, n)
\end{eqnarray*}
and
\begin{eqnarray*}
(c\cdot b)\cdot a & \cong & (t(tc+ub)+ua)(\mbox{mod}\, n)\\
 & \cong & (t^{2}c+tub+ua)(\mbox{mod}\, n)
\end{eqnarray*}

\noindent By using hypothesis, we get $(a\cdot b)\cdot c=(c\cdot b)\cdot a\,\,\forall\, a,b,c\in Z_{n}$.

\noindent Now we show that the class $(Z_{n},\cdot)$ is nonassociative
in general:
\begin{eqnarray*}
(a\cdot b)\cdot c & \cong & (t^{2}a+tub+uc)(\mbox{mod}\, n),
\end{eqnarray*}
\begin{eqnarray*}
a\cdot(b\cdot c) & \cong & (ta+utb+u^{2}c)(\mbox{mod}\, n).
\end{eqnarray*}

\noindent Since this is not necessary that
\begin{eqnarray*}
t^{2}a+tub+uc & \cong & (ta+utb+u^{2}c)(\mbox{mod}\, n).
\end{eqnarray*}

\noindent Hence $Z_{n}(t,u)$ is an AG-groupoid which may or may not
be associative.\end {proof}

\noindent We denote this AG-groupoid by $\{Z_{n},(t,u),\cdot\}$-AG-groupoid
$(\mbox{mod }n)$ or in short by $Z_{n}(t,u)$-AG-groupoid $(\mbox{mod }n)$.
For varying values of $t$ and $u$ and by putting some conditions
on $t$ and $u$, we get different classes of AG-groupoids $(\mbox{mod }n)$
for some fixed integer $n\geq3$. These new classes of AG-groupoids
$(\mbox{mod }n)$ will be denoted by $Z_{AG}^{*}(n)$, $Z_{AG}^{**}(n)$
and $Z_{AG}^{***}(n)$. In the following we list some examples to
show the existence of these modulo AG-groupoids.

\noindent \begin{enumerate}[(i)]

\item $Z_{3}(2,1)$ in $Z(3)$ is an AG-groupoid (in fact an AG-group):

\begin{center}
\begin{tabular}{l|lll}
$\cdot$  & $0$  & $1$  & $2$ \tabularnewline
\hline 
$0$  & $0$ & $1$  & $2$ \tabularnewline
$1$  & $2$  & $0$ & $1$\tabularnewline
$2$  & $1$  & $2$  & $0$ \tabularnewline
\end{tabular}
\par\end{center}

\item $Z_{8}(6,4)$ in $Z^{*}(8)$ is an AG-groupoid:

\begin{center}
\begin{tabular}{l|llllllll}
$\cdot$  & $0$  & $1$  & $2$  & $3$ & $4$ & $5$ & $6$ & $7$\tabularnewline
\hline 
$0$  & $0$ & $4$ & $0$ & $4$ & $0$ & $4$ & $0$ & $4$\tabularnewline
$1$  & $6$ & $2$ & $6$ & $2$ & $6$ & $2$ & $6$ & $2$\tabularnewline
$2$  & $4$ & $0$ & $4$ & $0$ & $4$ & $0$ & $4$ & $0$\tabularnewline
$3$ & $2$ & $6$ & $2$ & $6$ & $2$ & $6$ & $2$ & $6$\tabularnewline
$4$ & $0$ & $4$ & $0$ & $4$ & $0$ & $4$ & $0$ & $4$\tabularnewline
$5$ & $6$ & $2$ & $6$ & $2$ & $6$ & $2$ & $6$ & $2$\tabularnewline
$6$ & $4$ & $0$ & $4$ & $0$ & $4$ & $0$ & $4$ & $0$\tabularnewline
$7$ & $2$ & $6$ & $2$ & $6$ & $2$ & $6$ & $2$ & $6$\tabularnewline
\end{tabular}
\par\end{center}

\item$Z_{AG}^{**}(3)=\{Z_{3}(1,1),Z_{3}(2,1)\}$.

\item$Z_{AG}^{***}(4)=\{Z_{4}(1,1),Z_{4}(2,0),Z_{4}(3,1)\}$.

\item$Z_{AG}^{**}(5)=\{Z_{5}(1,1),Z_{5}(4,1),Z_{5}(3,4),Z_{5}(2,4)\}$. 

\item$Z_{AG}^{**}(6)=\{Z_{6}(1,1),Z_{6}(2,4),Z_{6}(3,3),Z_{6}(4,4),Z_{6}(5,1)\}$
and so on.

\end{enumerate}

\noindent We immediately have the following consequences of Theorem
\ref{GPDS}. 

\noindent \begin {corollary} Any AG-groupoid in $Z_{AG}^{***}(n)$
is a commutative semigroup if $t=u$.\end{corollary}

\noindent \begin {proof}If $t=u$ then the binary operation becomes
commutative which forces associativity in the AG-groupoid.\end {proof}
\begin{example}
\noindent \label{T3}$Z_{6}(4,4)$ in $Z_{AG}^{**}(6)$ is a commutative
semigroup given by the table:
\end{example}
\begin{center}
\begin{tabular}{l|llllll}
$\cdot$  & $0$  & $1$  & $2$  & $3$ & $4$ & $5$\tabularnewline
\hline 
$0$  & $0$ & $4$ & $2$ & $0$ & $4$ & $2$\tabularnewline
$1$  & $4$ & $2$ & $0$ & $4$ & $2$ & $0$\tabularnewline
$2$  & $2$ & $0$ & $4$ & $2$ & $0$ & $4$\tabularnewline
$3$ & $0$ & $4$ & $2$ & $0$ & $4$ & $2$\tabularnewline
$4$ & $4$ & $2$ & $0$ & $4$ & $2$ & $0$\tabularnewline
$5$ & $2$ & $0$ & $4$ & $2$ & $0$ & $4$\tabularnewline
\end{tabular}
\par\end{center}

\noindent Next we characterize AG-groupoids (mod $n$).

\noindent \begin {theorem} \label{the2} An AG-groupoid in  $Z_{AG}^{**}(n)$
is a $T^{3}$-AG-groupoid, if $t=u$.\end {theorem}

\noindent \begin {proof} Let $t=u$, then to show that an AG-groupoid
in $Z_{AG}^{**}(n)$ is a $T^{3}$-AG-groupoid, we will have to show
that it is a $T_{l}^{3}$-AG-groupoid as well as a $T_{r}^{3}$-AG-groupoid. 

\noindent For $T_{l}^{3}$-AG-groupoid, let 
\begin{eqnarray*}
a\cdot b & = & a\cdot c\\
\Rightarrow(ta+ub) & \cong & (ta+uc)(\mbox{mod}\, n)\\
\Rightarrow ub & \cong & uc(\mbox{mod }n)\\
\Rightarrow tb & \cong & tc(\mbox{mod }n)\qquad(2.1)
\end{eqnarray*}
Now 
\begin{eqnarray*}
b\cdot a & \cong & (tb+ua)(\mbox{mod}\, n)\\
 & \cong & (tc+ua)(\mbox{mod}\, n)\qquad(\mbox{by Equation 2.1})\\
\Rightarrow b\cdot a & = & c\cdot a
\end{eqnarray*}
Hence an AG-groupoid in $Z_{AG}^{**}(n)$ is $T_{l}^{3}$-AG-groupoid.
Similarly we can show that an AG-groupoid in $Z_{AG}^{**}(n)$ is
$T_{r}^{3}$-AG-groupoid. Hence any AG-groupoid in $Z_{AG}^{**}(n)$
is $T^{3}$-AG-groupoid if $t=u$.\end {proof}

\noindent \begin {theorem} Every AG-groupoid in $Z_{AG}^{*}(n)$
is a $T^{3}$-AG-groupoid, if $n$ is prime.\end {theorem}

\noindent \begin {proof} Let $n$ be any prime number then to show
that an AG-groupoid in $Z_{AG}^{*}(n)$ is a $T^{3}$-AG-groupoid,
we will have to show that it is a $T_{l}^{3}$-AG-groupoid as well
as a $T_{r}^{3}$-AG-groupoid.

\noindent For $T_{l}^{3}$-AG-groupoid, let $a,b,c\in G$, and 
\begin{eqnarray*}
a\cdot b & = & a\cdot c\\
\Rightarrow(ta+ub) & \cong & (ta+uc)(\mbox{mod}\, n)\\
\Rightarrow ub & \cong & uc(\mbox{mod }n)\\
\Rightarrow u(b-c) & \cong & 0(\mbox{mod }n)
\end{eqnarray*}

\noindent as $n\nmid u$, because $n$ is a prime number. Therefore,
$n\mid(b-c)\Rightarrow b\cong c(\mbox{mod }n)$, and consequently;
\begin{eqnarray*}
b\cdot a & \cong & (tb+ua)(\mbox{mod}\, n)\\
 & \cong & (tc+ua)(\mbox{mod}\, n)\\
\Rightarrow b\cdot a & = & c\cdot a
\end{eqnarray*}
Hence every AG-groupoid $G$ in $Z_{AG}^{*}(n)$ is $T_{l}^{3}$-AG-groupoid.
Similarly we can show that every AG-groupoid in $Z_{AG}^{*}(n)$ is
$T_{r}^{3}$-AG-groupoid. Hence every AG-groupoid in $Z_{AG}^{*}(n)$
is $T^{3}$-AG-groupoid if $n$ is prime.\end {proof}

\noindent \begin {example} $Z_{5}(3,4)$ in $Z_{AG}^{*}(5)$ is a
$T^{3}$-AG-groupoid:\end {example}

\begin{center}
\begin{tabular}{l|lllll}
$\cdot$  & $0$  & $1$  & $2$  & $3$ & $4$\tabularnewline
\hline 
$0$  & $0$ & $4$ & $3$ & $2$ & $1$\tabularnewline
$1$  & $3$ & $2$ & $1$ & $0$ & $4$\tabularnewline
$2$  & $1$ & $0$ & $4$ & $3$ & $2$\tabularnewline
$3$ & $4$ & $3$ & $2$ & $1$ & $0$\tabularnewline
$4$ & $2$ & $1$ & $0$ & $4$ & $3$\tabularnewline
\end{tabular}
\par\end{center}

\noindent Also in Example \ref{T3}; $Z_{6}(4,4)$ in $Z_{AG}^{*}(6)$
is a $T^{3}$-AG-groupoid. However, the result is not true in general.
For example, $Z_{8}(6,4)$ is not a $T^{3}$-AG-groupoid:

\begin{center}
\begin{tabular}{l|llllllll}
$\cdot$  & $0$ & $1$ & $2$ & $3$ & $4$ & $5$ & $6$ & $7$\tabularnewline
\hline 
$0$ & $0$ & $4$ & $0$ & $4$ & $0$ & $4$ & $0$ & $4$\tabularnewline
$1$ & $6$ & $2$ & $6$ & $2$ & $6$ & $2$ & $6$ & $2$\tabularnewline
$2$ & $4$ & $0$ & $4$ & $0$ & $4$ & $0$ & $4$ & $0$\tabularnewline
$3$ & $2$ & $6$ & $2$ & $6$ & $2$ & $6$ & $2$ & $6$\tabularnewline
$4$ & $0$ & $4$ & $0$ & $4$ & $0$ & $4$ & $0$ & $4$\tabularnewline
$5$ & $6$ & $2$ & $6$ & $2$ & $6$ & $2$ & $6$ & $2$\tabularnewline
$6$ & $4$ & $0$ & $4$ & $0$ & $4$ & $0$ & $4$ & $0$\tabularnewline
$7$ & $2$ & $6$ & $2$ & $6$ & $2$ & $6$ & $2$ & $6$\tabularnewline
\end{tabular}
\par\end{center}

The following theorem shows that $Z_{AG}^{*}(n)$ is a subclass of
transitively commutative AG-groupoid.

\noindent \begin {theorem} \label{th3} Every AG-groupoid in $Z_{AG}^{*}(n)$
is transitively commutative AG-groupoid.\end {theorem}

\noindent \begin {proof} To show that every AG-groupoids in $Z_{AG}^{*}(n)$,
is transitively commutative AG-groupoid it is sufficient if we show
that an arbitrary AG-groupoid $H$ in $Z_{AG}^{*}(n)$ is transitively
commutative AG-groupoid. Now for any $a,b,c\in H,$ we show that for
$a\cdot b=b\cdot a\,\,\mbox{and}\,\, b\cdot c=c\cdot b\Rightarrow a\cdot c=c\cdot a.$
Let 
\begin{eqnarray*}
a\cdot b & = & b\cdot a\\
\Rightarrow(ta+ub) & \cong & (tb+ua)(\mbox{mod}\, n)\\
\Rightarrow t(a-b)+u(b-a) & \cong & 0(\mbox{mod }n)
\end{eqnarray*}

\noindent Similarly,
\begin{eqnarray*}
b\cdot c & = & c\cdot b\\
\Rightarrow(tb+uc) & \cong & (tc+ub)(\mbox{mod}\, n)\\
\Rightarrow t(b-c)+u(c-b) & \cong & 0(\mbox{mod }n)
\end{eqnarray*}

\noindent as $n\mid t(a-b)+u(b-a)$ and $n\mid t(b-c)+u(c-b)\Rightarrow n\mid t(a-b)+u(b-a)+t(b-c)+u(c-b)$
\begin{eqnarray*}
\Rightarrow t(a-b)+u(b-a)+t(b-c)+u(c-b) & \cong & 0(\mbox{mod}\, n)\\
t(a-c)+u(-a+c) & \cong & 0(\mbox{mod}\, n)\\
(ta+uc)-(tc+ua) & \cong & 0(\mbox{mod}\, n)\\
ta+uc & \cong & (tc+ua)(\mbox{mod}\, n)\\
\Rightarrow ac & = & ca.
\end{eqnarray*}

\noindent Hence every AG-groupoids in $Z_{AG}^{*}(n)$ is transitively
commutative AG-groupoid.\end {proof}

\noindent \begin {example} $Z_{7}(5,4)$ is transitively commutative
AG-groupoid:\end {example}

\begin{center}
\begin{tabular}{l|lllllll}
$\cdot$  & $0$ & $1$ & $2$ & $3$ & $4$ & $5$ & $6$\tabularnewline
\hline 
$0$ & $0$ & $4$ & $1$ & $5$ & $2$ & $6$ & $3$\tabularnewline
$1$ & $5$ & $2$ & $6$ & $3$ & $0$ & $4$ & $1$\tabularnewline
$2$ & $3$ & $0$ & $4$ & $1$ & $5$ & $2$ & $6$\tabularnewline
$3$ & $1$ & $5$ & $2$ & $6$ & $3$ & $0$ & $4$\tabularnewline
$4$ & $6$ & $3$ & $0$ & $4$ & $1$ & $5$ & $2$\tabularnewline
$5$ & $4$ & $1$ & $5$ & $2$ & $6$ & $3$ & $0$\tabularnewline
$6$ & $2$ & $6$ & $3$ & $0$ & $4$ & $1$ & $5$\tabularnewline
\end{tabular}
\par\end{center}

\noindent \begin {theorem} Every AG-groupoid in $Z_{AG}^{*}(n)$
is a cancellative AG-groupoid, if $n$ is prime.\end {theorem}

\noindent \begin {proof} To show that for any prime number $n$;
every AG-groupoid in $Z_{AG}^{*}(n)$ is a cancellative AG-groupoid,
it is sufficient if we show an arbitrary AG-groupoid is left cancellative
AG-groupoid. 

\noindent For left cancellative AG-groupoid, let 
\begin{eqnarray*}
a\cdot x & = & a\cdot y\\
\Rightarrow(ta+ux) & \cong & (ta+uy)(\mbox{mod}\, n)\\
\Rightarrow u(x-y) & \cong & 0(\mbox{mod }n)
\end{eqnarray*}
as $n\nmid u$, because $n$ is a prime number. Therefore, it means
that $n\mid(x-y)\Rightarrow x\cong y(\mbox{mod }n)$. Hence every
AG-groupoid is left cancellative. As every left cancellative AG-groupoid
is right cancellative AG-groupoid \cite{MshahT}. Hence every AG-groupoid
in $Z_{AG}^{*}(n)$ is a cancellative AG-groupoid.\end {proof}

\noindent \begin {example} $Z_{5}(3,4)$ in $Z_{AG}^{*}(5)$ is a
cancellative AG-groupoid:\end {example}

\begin{center}
\begin{tabular}{l|lllll}
$\cdot$  & $0$ & $1$ & $2$ & $3$ & $4$\tabularnewline
\hline 
$0$ & $0$ & $4$ & $3$ & $2$ & $1$\tabularnewline
$1$ & $3$ & $2$ & $1$ & $0$ & $4$\tabularnewline
$2$ & $1$ & $0$ & $4$ & $3$ & $2$\tabularnewline
$3$ & $4$ & $3$ & $2$ & $1$ & $0$\tabularnewline
$4$ & $2$ & $1$ & $0$ & $4$ & $3$\tabularnewline
\end{tabular}
\par\end{center}

\noindent \begin {theorem} An AG-groupoid in $Z_{AG}^{***}(n)$ is
an AG-band if $t+u=1$. \end {theorem}

\noindent \begin {proof}Let $t+u=1$, to show that $Z_{AG}^{***}(n)$
is an AG-band it is sufficient to show that $a\cdot a\cong a(\mbox{mod }n)$; 

\noindent 
\begin{eqnarray*}
a\cdot a & \cong & (ta+ua)(\mbox{mod \ensuremath{n)}}\\
 & \cong & a(t+u)(\mbox{mod \ensuremath{n})}\\
 & \cong & a(\mbox{mod \ensuremath{n})}
\end{eqnarray*}
Hence the claim.\end {proof}

\noindent \begin {example}$Z_{5}(2,4)$ is an AG-band:\end {example}

\begin{center}
\begin{tabular}{l|lllll}
$\cdot$  & $0$  & $1$  & $2$  & $3$ & $4$\tabularnewline
\hline 
$0$  & $0$ & $4$ & $3$ & $2$ & $1$\tabularnewline
$1$  & $2$ & $1$ & $0$ & $4$ & $3$\tabularnewline
$2$  & $4$ & $3$ & $2$ & $1$ & $0$\tabularnewline
$3$ & $1$ & $0$ & $4$ & $3$ & $2$\tabularnewline
$4$ & $3$ & $2$ & $1$ & $0$ & $4$\tabularnewline
\end{tabular}
\par\end{center}

\section{\noindent Existence of AG-Groups $(\mbox{mod}\, n)$}

In this section, we introduce a special class of AG-groupoids $(\mbox{mod }n)$,
namely AG-groups $(\mbox{mod }n)$ and give some of its characterizations.
The following theorem shows the existence of AG-groups $(\mbox{mod }n)$
for $n\geq3$, and indeed it gives a simple way of construction of
AG-groups of any finite order.

\noindent \begin {theorem} \label{AGd} Let $Z_{n}=\{0,1,2,\cdots,n-1\}\, n\geq3,\, n<\infty$.
A groupoid in $Z^{**}(n)$ is an AG-group if $t^{2}\cong1(\mbox{mod}\, n)\,\,\mbox{for}\, t\in Z_{n}\setminus\{0\}$.\end{theorem}

\noindent \begin {proof}Given that $Z_{n}=\{0,1,2,\cdots,n-1\}\, n\geq3,\, n<\infty$;
satisfies $t^{2}\cong1(\mbox{mod}\, n)$ for $t\in Z_{n}\setminus\{0\}$,
we have to show $Z_{n}(t,1)$ in $Z^{**}(n)$ is an AG-group.

\noindent We show that the left invertive law, $(a\cdot b)\cdot c=(c\cdot b)\cdot a$
holds in $Z_{n}$; 

\noindent 
\begin{eqnarray*}
(a\cdot b)\cdot c & \cong & (t(ta+b)+c)(\mbox{mod}\, n)\\
 & \cong & (t^{2}a+tb+c)(\mbox{mod}\, n)
\end{eqnarray*}
and
\begin{eqnarray*}
(c\cdot b)\cdot a & \cong & (t(tc+b)+a)(\mbox{mod}\, n)\\
 & \cong & (t^{2}c+tb+a)(\mbox{mod}\, n)
\end{eqnarray*}

\noindent Hence $Z_{n}(t,1)$ is an AG-groupoid as $t^{2}\cong1(\mbox{mod}\, n)$
and $(a\cdot b)\cdot c=(c\cdot b)\cdot a$.

\noindent \emph{Existence of left identity:} `$0$' is the correspondent
left identity;
\begin{eqnarray*}
0\cdot x & \cong & x(\mbox{mod}\, n)x\,\forall\, x\in Z_{n}
\end{eqnarray*}
but~
\begin{eqnarray*}
x\cdot0 & \cong & (tx)(\mbox{mod}\, n).
\end{eqnarray*}

\noindent \emph{Existence of inverses:} $(n-1)tx=-tx$ is the inverse
of $x\,\forall\, x\in Z_{n}$;
\begin{eqnarray*}
(-tx)\cdot x & \cong & (t(-tx)+x)(\mbox{mod}\, n)\\
 & \cong & (-(t^{2}-1)x)(\mbox{mod}\, n)\\
 & \cong & 0(\mbox{mod}\, n)
\end{eqnarray*}

\noindent and
\begin{eqnarray*}
x\cdot(-tx & )\cong & (tx+(-tx))(\mbox{mod}\, n)\\
 & \cong & 0(\mbox{mod}\, n)
\end{eqnarray*}

\noindent Hence $(Z_{n},\cdot)$ is an AG-group (mod $n$).\end {proof}

\noindent We denote this AG-group $(\mbox{mod }n)$ by $\{Z_{n},(t,1),\cdot\}\mbox{-AG-group}\,(\mbox{mod }n)$
or in short by $Z_{n}(t,1)$-AG-group $(\mbox{mod }n)$. 

\noindent \begin {corollary} Any AG-groupoid in $Z_{AG}^{**}(n)$
is an abelian group if $t=1$.\end{corollary}

\noindent \begin {proof} If $t=1$ then 1 becomes the identity of
the $Z_{n}(1,1)$-AG-group $(\mbox{mod }n)$ and so it becomes abelian
group by \cite[Theorem 2]{MshahT}. \end {proof}

\noindent \begin {corollary} Let $Z_{n}=\{0,1,2,\cdots,n-1\},\, n\geq3,\, n<\infty$.
Then $Z_{n}(n-1,1)$ is an AG-group $(\mbox{mod }n)$.\end{corollary}

\noindent \begin {proof} Since $(n-1)^{2}\cong1(\mbox{mod}\, n)$.
The proof now follows by Theorem \ref{AGd}.\end {proof}

\noindent We denote this AG-group by $\{Z_{n},(n-1,1),\cdot\}$-AG-group
$(\mbox{mod }n)$ or in short by $Z_{n}(n-1,1)$-AG-group $(\mbox{mod }n)$.

\end{document}